\documentclass[11pt,êeqno]{article}
\usepackage{amsfonts}
\usepackage[tbtags]{amsmath}
\topskip=\baselineskip%
\newcommand{\nty}{n \to \infty}

\newcommand{\lr}{\left(}
\newcommand{\lb}{\left[}
\newcommand{\rb}{\right]}
\newcommand{\rr}{\right)}

\newcommand{\lv}{\left|}
\newcommand{\rv}{\right|}

\newcommand{\rp}{\right.}
\newcommand{\kn}{k_n}

\newcommand{\eel}{\end{lemma}}

\newcommand{\xin}{X_{i:n}}
\newcommand{\uin}{U_{i:n}}
\newcommand{\xia}{\xi_p}
\newcommand{\qa}{\xi_{p_n}}
\newcommand{\qan}{\xi_{p_n n:n}}
\newcommand{\an}{p_n}
\newcommand{\xian}{\xi_{p\, n:n}}
\newcommand{\xiann}{\xi_{p_n\, n:n}}

\newcommand{\nau}{N_{p,u }}
\newcommand{\lon}{\log n/n}

\newcommand{\alp}{p}

\newcommand{\aln}{p_n}

\newcommand{\na}{N_{p}}

\numberwithin{equation}{section}

\newtheorem{theorem}{Theorem}[section]

\newtheorem{lemma}{Lemma}[section]

\newtheorem{corollary}{Corollary}[section]

\newtheorem{example}{Example}[section]
\newtheorem{remark}{Remark}[section]

\title{{\huge \sf On the Bahadur -- Kiefer Representation for \\  Intermediate Sample Quantiles}
}
\author{{\sf Nadezhda Gribkova$^{*1}$} \ and  \ \ {\sf Roelof Helmers$^{2}$}}
\date{
{\small \it
\centerline{$^{1}$~St.Petersburg~State~University,~Mathematics~and~Mechanics~Faculty,}
\centerline{198504,~St.Petersburg,~Stary~Peterhof,~Universitetsky~pr.~28,~Russia}
\centerline{E-mail:~nv.gribkova@gmail.com}
\centerline{$^{2}$ Center~for~Mathematics~and~Computer~Science}
\centerline{P.O.Box~94079,~1090~GB~Amsterdam,~The~Netherlands}
\centerline{E-mail:~helmers@cwi.nl}
}} 
\begin{document}
\maketitle 
       
\begin{abstract}
{\small \it  We investigate a~Bahadur-Kiefer type representation for the $\aln$-th empirical quantile corresponding to a~sample of $n$ i.i.d. random variables, when  $\aln \in (0,1)$ is a sequence which, in particular, may tend to $0$ or $1$, i.e. we consider the  case of intermediate sample quantiles. We obtain an~'in probability' version of the Bahadur -- Kiefer type representation for a~$\kn$-th order statistic when $r_n=\kn\wedge (n-k_n)\to \infty$ under some mild regularity conditions, and an 'almost sure' version under additional assumption  that $\log n/r_n\to 0$, $\nty$.
A~representation for the~sum of order statistics laying between  the population $\aln$-quantile and the corresponding empirical quantile is also established.}
\end{abstract}
\begin{quote}
\end{quote}

\section{Introduction}
The classical Bahadur -- Kiefer  representation was established by Bahadur \cite{b} and Kiefer  \cite{kiefer67}-\cite{kiefer70a}, it allows one to replace the quantile process by (-1) times the empirical process with an almost sure uniform error of the order $n^{-1/4+o(1)}$, where $n$ is the real i.i.d. data sample size (see, e.g.,  Shorack and Wellner \cite{shorack}, Deheuvels and Mason \cite{deh_mas90}, Deheuvels \cite{deh_00}, see also references therein).

In this paper we investigate the asymptotic behavior of the so-called intermediate sample quantile, i.e. of the $k_n$-th order statistic, $1\le k_n\le n$, when $r_n:\,=k_n\wedge (n-k_n)\to\infty$, $p_n:\,=k_n/n\to 0$ (or $p_n\to 1$), as $\nty$. We obtain Bahadur--Kiefer type representations for intermediate sample  quantiles under a~mild regularity condition, and we establish also a~representation for sum of the order statistics laying between the population $p_n$-th quantile and the corresponding sample quantile.

Consider a sequence $X_1,X_2,\dots $ of independent identically
distributed (i.i.d.) real-valued random variables
(r.v.) with common distribution function ($df$)~$F$, and for each
integer $n \ge 1$ let \ $X_{1:n}\le \dots \le X_{n:n}$ denote the
order statistics based on the sample $X_1,\dots ,X_n$. Let $F^{-1}(u)= \inf \{ x: F(x) \ge u \}$, \ $0<u\le 1$, \
$F^{-1}(0)=F^{-1}(0^+)$, denote the left-continuous inverse function of $df$ $F$,  and  $F_n$,  $F_n^{-1}$ ---  the
empirical $df$ and its inverse respectively, put $f=F'$ to be a density of the underlying distribution when it exists.
 Let  $\xia=F^{-1}(\alp)$, \  $\xian=F_n^{-1}(\alp)$ \ denote $\alp$-th  population  and  sample  quantile respectively.

 For a~fixed  $\alp\in (0,1)$ assuming that $F$ has at least two continuous derivatives in a~neighborhood of $\xi_p$ and $f(\xia)>0$, Bahadur \cite{b} first establish the almost sure result:
\begin{equation}
\label{bahadur_b}
\xian=\xia-\frac{F_n(\xia) -\alp}{f(\xia)}+R_n(p),
\end{equation}
where  $R_n(p) =O_{\text{a.s.}}\lr n^{-3/4}(\log n)^{1/2}(\log\log n)^{1/4}\rr$ (a sequence of random variables $R_n$ is said to be $O_{\text{a.s.}}(\tau_n)$ if $R_n/\tau_n$  is almost surely bounded). Kiefer in a~sequence of papers \cite{kiefer67}-\cite{kiefer70a} proved that if $f'$ is bounded in
a~neighborhood of $\alp$ and $f(\xia)>0$, then
$
\limsup_{\nty}\pm\, n^{3/4} (\log\log n)^{-3/4}R_n(p)=\frac{2^{5/4}3^{-3/4}(\alp (1-\alp))^{1/4}}{f(\xia)}
$
$a.s.$ for either choice of sign. In Reiss \cite{r} a~version of Bahadur's result with a~remainder term, which is of the order $O(\lon)^{3/4}$ in probability was obtained: if the density $f=F'$ is Lipschitz in a~neighborhood of $\alp$ and $f(\xia)>0$, then \eqref{bahadur_b}  holds true and  $\boldsymbol{P}\lr |R_n(p)|>A(\lon)^{3/4} \rr\le Bn^{-c}$ for every $c>0$, where  $A,B$ are some positive constants, not depending on $n$.

Our interest in Bahadur-Kiefer type representation for intermediate empirical quantile was first motivated by its uses in the second order asymptotic analysis of trimmed sums. It turns  out (see Gribkova and Helmers \cite{gh2006}-\cite{gh2010})  that the Bahadur -- Kiefer properties provide a~very useful tool in investigation of the asymptotic behavior of the distributions of  trimmed sums of i.i.d. r.v.'s, slightly trimmed sums and their studentized versions. In particular, the Bahadur's representation allows us to construct a~$U$-statistics type stochastic approximation for these statistics, which will enable us to establish  the Berry -- Esseen type bounds and the Edgeworth expansions in Central Limit Theorems for normalized and studentized slightly trimmed sums.

We would like to emphasize, that the Bahadur-Kiefer type representation we obtain for a sum of order statistics lying between  the $\an$-th population quantile and the corresponding empirical quantile  (cf. Theorem~\ref{stm_Bah_2}), is  especially useful in the construction of the $U$-statistic type approximation for a~(slightly) trimmed sum, as it provides a~quadratic term of the desired $U$-statistic. Note also that formally the representation~\eqref{stm_(A.2)} (cf. Theorem~\ref{stm_Bah_2}) can be obtained by integrating of the corresponding Bahadur -- Kiefer process in interval $[\xi_{\aln n: n},\xi_{\aln})$, however we prove representation  \eqref{stm_(A.2)} for  intermediate order statistics (i.e. when $p_n\to 0$ (or $p_n\to 1$). The remainder terms in our representations are shown to be of a~suitable order of magnitude similar as in Reiss \cite{r}.

Part of our results can be compared with an~earlier result obtained  by Chanda \cite{Chanda_92}, who established  the Bahadur -- Kiefer representation for the intermediate $k_n$-th order statistics, assuming the somewhat restrictive condition  $n^a/k_n\to 0$ for some $a>0$ and, in addition, some strong regularity conditions on $F$ must be satisfied.

We conclude this introduction by noting that some extensions of Bahadur's result to dependent random variables have been proved by Sen \cite{sen68} (cf.~also Wu \cite{Wu}). The validity of Bahadur's representation for a~bootstrapped $\alp$-quantile was proved (as an~auxiliary result) in Gribkova and Helmers \cite{gh2007}.  Deheuvels \cite{deh_09} established a~multivariate Bahadur--Kiefer representation for the empirical copula process.
\section{Statement of results}
Let $\kn$  be a sequences of integers, such that $0\le  \kn \le n$, and $r_n=\kn\wedge (1-k_n) \to \infty$, as $\nty$. Put
$\aln=\kn/n$, and let $\qa=F^{-1}(\aln)$, \  $\xiann=F_n^{-1}(\aln)$ \ denote $\aln$-th  population  and  empirical quantile respectively.

Define two  numbers
\begin{equation}
\label{seq}
0\le a_1= \liminf_{\nty}\, \an \ \ \le \ \ a_2= \limsup_{\nty}\, \an  \le 1. \\
\end{equation}
We will assume throughout this note that the following
smoothness condition is satisfied.

\medskip \noindent
$[A_1].$ \ \ {\it The function  $F^{-1}$ is differentiable in some open set \ $U\subset (0,1)$ (i.e. the density $f=F'$ exists and is positive in  $F^{-1}(U)$), moreover
\begin{equation}
\label{1.6}
\begin{array}{cccccccccc}
 & (0,\varepsilon), \ \ \ \text{if }\  0=a_1=a_2, & &(1-\varepsilon,1),  \ \ \text{if} \ \ a_1=a_2=1,\\
U  \supset & (0,a_2]\ \ \ \text{if }\ \ 0=a_1<a_2, &\qquad U \supset &\ \ \ \ [a_1,1), \ \ \ \ \ \text{if} \ \ 0<a_1<a_2=1,\\
& \ \ \ \ \   [a_1,a_2], \  \text{if}\ \   0<a_1\le a_2<1, & & \ \ \ (0,1),  \ \ \ \ \ \ \text{if} \ \ a_1=0,\ a_2=1,\\
\end{array}
\end{equation}
with some $0<\varepsilon \le 1$ in cases given in the first lines of \eqref{1.6}).}

To state our results we will need also the following condition:

\medskip \noindent
$[A_2].$ \qquad \qquad \qquad \qquad $\qquad r_n^{-1}\log n\to 0, \quad \ \ \nty.$

\medskip
Let $h$ be a real-valued function defined on the set $F^{-1}(U)$
(cf. \eqref{1.6}). Take an arbitrary $0<C<\infty$ and for all sufficiently large $n$ define
\begin{equation}
\label{1.8}
\begin{split}
\Psi_{\an,h}(C)&=\sup_{|t|\le C}\lv h\circ F^{-1}\Bigl( \an
+t\sqrt{\frac{r_n \log r_n}{n^2}}\Bigr) -h\circ F^{-1}\Bigl(\an\Bigr)\rv,
\end{split}
\end{equation}
where $h\circ F^{-1}(u)=h\lr F^{-1}(u)\rr$. Note that $\an
+t\sqrt{\frac{r_n \log r_n}{n^2}}=\an  +t\, \frac{r_n}n \, \sqrt{\frac{\log r_n}{r_n}}=\an  +t\, \frac{r_n}n \,o(1)$,  $\nty$. In particular, this implies  that the function introduced in~\eqref{1.8} is
well-defined for all sufficiently large~$n$.

Next we define a function $\widehat{\Psi}_{\an,h}(C)$ which is equal to $\Psi_{\an,h}(C)$, where $\log r_n$ is replaced by $\log n$. Similarly as before we show that it is well-defined for all sufficiently large~$n$ if condition $[A_2]$ holds true.

   We will obtain the Bahadur-Kiefer type representations for some smooth function of the empirical quantile,
as it turned out (cf.~\cite{gh2006}-\cite{gh2010}) that these extensions are very useful in construction of the $U$-statistic type stochastic approximations for the trimmed sums.

  Let $G(x)$, $x\in R$, be a real-valued function, $g=G'$ -- its derivative when it exists,
and let $(g/f)(x)$ and $(|g|/f)(x)$ denote the ratios $g(x)/f(x)$ and $|g(x)|/f(x)$ respectively.
\begin{theorem}
\label{stm_Bah_1}
Suppose that $r_n\to \infty$, as $\nty$, the condition $[A_1]$ holds true and $G$ is
differentiable on the set $F^{-1}(U)$.
Then
\begin{equation}
\label{stm_(A.1)}
G(\qan)-G(\qa)=-[F_n(\qa)-F(\qa)]\frac gf(\qa)+R_n(p_n),
\end{equation}
where  for each $c>0$
\begin{equation}
\label{delta_(1.1)}
\boldsymbol{P}(|R_n(p_n)|>\Delta_n)=O\lr r_n^{-c}\rr,
\end{equation}
 with
\begin{equation*}
\Delta_n=A\, (\an(1-\an))^{1/4}\lr\frac{\log r_n}n\rr^{3/4}\frac {|g|}f(\qa)+ B\,(\an(1-\an))^{1/2}\lr\frac{\log r_n}n\rr^{1/2}   \Psi_{\an,\frac gf}\, (C),
\end{equation*}
where $A$, $B$ and $C$ are some positive constants, which  depend only
on $c$.

Moreover, if additionally the condition $[A_2]$ is also satisfied, then \eqref{stm_(A.1)} holds true and
\begin{equation}
\label{delta_h_(1.1)}
\boldsymbol{P}(|R_n(p_n)|>\widehat{\Delta}_n)=O\lr n^{-c}\rr,
\end{equation}
for each $c>0$ with
\begin{equation*}
\widehat{\Delta}_n=A\, (\an(1-\an))^{1/4}\lr\frac{\log n}n\rr^{3/4}\frac {|g|}f(\qa)+ B\,(\an(1-\an))^{1/2}\lr\frac{\log n}n\rr^{1/2}   \widehat{\Psi}_{\an,\frac gf}\, (C),
\end{equation*}
where $A$, $B$ and $C$ are some positive constants, which  depend only
on $c$.
\end{theorem}

Theorem~\ref{stm_Bah_1} is a Bahadur-Kiefer type result. For the~special case when $0<\alp<1$ is fixed it is stated in Lemma~3.1 of~\cite{gh2006} (cf.~also Lemma~4.1,~\cite{gh2007} and  Reiss~\cite{r}).
\begin{remark} \label{rem_1} It is easy to see that if one compares the first term on the r.h.s. of \eqref{stm_(A.1)}  and the orders of magnitude of the quantities $\Delta_n$, \ $\widehat{\Delta}_n$ given in \eqref{delta_(1.1)}---\eqref{delta_h_(1.1)}  that  relation \eqref{stm_(A.1)} provides a~representation with a~remainder term $R_n(p_n)$ of smaller order than the first term if and only if $ \Psi_{\an,\frac gf}\, (C)=o\bigl(\frac {|g|}f(\qa)\bigr)$ and $ \widehat{\Psi}_{\an,\frac gf}\, (C)=o\bigl(\frac {|g|}f(\qa)\bigr)$ for every fixed $C>0$, as \ $\nty$. The same remark is valid for the~two assertions stated in Theorem~\ref{stm_Bah_2} below.
\end{remark}

We relegate  proofs of our results  to  sections~\ref{proof1}~--~\ref{proof2}.
\begin{theorem}
\label{stm_Bah_2}   Suppose that $r_n\to \infty$, as $\nty$, the condition $[A_1]$ holds true and $G$ is
differentiable on the set $F^{-1}(U)$.
Then
\begin{equation}
\label{stm_(A.2)}
\int_{\qan}^{\qa}(G(x)-G(\qa))\, d\, F_n(x)=-\frac
12[F_n(\qa)-F(\qa)]^2\frac gf(\qa)+R_n(p_n),
\end{equation}
where
\begin{equation}
\label{delta_(1.2)}
\boldsymbol{P}(|R_n(p_n)|>\Delta_n)=O\lr r_n^{-c}\rr,
\end{equation}
for each $c>0$ with
\begin{equation*}
\Delta_n=A\, (\an(1-\an))^{3/4}\lr\frac{\log r_n}n\rr^{5/4}\frac {|g|}f(\qa)+ B\,\an(1-\an)\frac{\log r_n}n  \Psi_{\an,\frac gf}\, (C),
\end{equation*}
where $A$, $B$ and $C$ are some positive constants, which  depend only
on $c$.

Moreover, if additionally the condition $[A_2]$ is also satisfied, then \eqref{stm_(A.1)} holds true, and
\begin{equation}
\label{delta_h_(1.2)}
\boldsymbol{P}(|R_n(p_n)|>\widehat{\Delta}_n)=O\lr n^{-c}\rr
\end{equation}
 for each $c>0$ with
\begin{equation*}
\widehat{\Delta}_n=A\, (\an(1-\an))^{3/4}\lr\frac{\log n}n\rr^{5/4}\frac {|g|}f(\qa)+ B\,\an(1-\an)\frac{\log n}n  \widehat{\Psi}_{\an,\frac gf}\, (C),
\end{equation*}
where $A$, $B$ and $C$ are some positive constants, which  depend only
on $c$.
\end{theorem}

Theorem~\ref{stm_Bah_2} extends  Lemma~4.3 from \cite{gh2007} (cf. also Lemma~3.2,~\cite{gh2006}), where it was proved for
a~fixed $\alp$ to the~case that  $\an$ is a~sequence which may tend to $0$ or to $1$. Note also that if both  conditions $[A_1]$ and $[A_2]$ are satisfied, then Theorems~\ref{stm_Bah_1}---\ref{stm_Bah_2} and an application of the Borel-Cantelly lemma imply an~almost sure result, i.e.  $R_n(p_n)=O_{a.s.}(\widehat{\Delta}_n)$, as $\nty$.

Next we will state some consequences of the Theorems~\ref{stm_Bah_1}---\ref{stm_Bah_2} where  the~remainder terms are  given in simpler form. Our first two consequences concern the Bahadur-Kiefer type representations for the central (not intermediate) order statistics.
\begin{corollary}
\label{cor1} Suppose that $0<a_1\le a_2<1$, the condition $[A_1]$ holds true and the functions $f=F'$ and $g=G'$ satisfies a~H\"{o}lder condition of the order $\alpha\ge 1/2$ on the set $F^{-1}(U)$. Then \eqref{stm_(A.1)} is valid and $\boldsymbol{P}(|R_n(p_n)|>A(\lon)^{3/4})=O\lr n^{-c}\rr$ for each $c>0$, where $A>0$ is some constant, not depending on $n$.
\end{corollary}
\begin{corollary}  \label{cor2} Suppose that the conditions of the Corollary~\ref{cor1} are satisfied. Then \eqref{stm_(A.2)} is valid and $\boldsymbol{P}(|R_n(p_n)|>A(\lon)^{5/4})=O\lr n^{-c}\rr$ for each $c>0$, where $A>0$ is some constant, not depending on $n$.
\end{corollary}
To prove Corollaries~\ref{cor1}---\ref{cor2} it suffices to note that the condition $0<a_1\le a_2<1$ implies that  $[A_2]$ is automatically satisfied, moreover, due to condition $[A_1]$ the density $f$ is bounded away from zero on the set $F^{-1}([a_1-\delta, \, a_2+\delta])$ with some $\delta>0$, and hence, the ratio $g/f$ satisfies a~H\"{o}lder condition of the order $\alpha\ge 1/2$ on this set.  Then an~application of H\"{o}lder's condition to the function $\Psi_{\an,\frac gf}\, (C)$ (cf.~\eqref{1.8}) proves both corollaries.

Next we state several  corollaries  for the intermediate sample  quantiles provided some regularity conditions are satisfied.

Note that the second terms of $\Delta_n$ and $\widehat{\Delta}_n$ in \eqref{delta_(1.1)}-\eqref{delta_h_(1.1)} and in \eqref{delta_(1.2)}-\eqref{delta_h_(1.2)}, involving the functions $\Psi_{\an,\frac gf}\, (C)$ and $\widehat{\Psi}_{\an,\frac gf}\, (C)$, depend on the asymptotic  properties of the ratio $g/f$, and  we can describe some sets of conditions allowing to absorb these second terms in the first ones. We will need the following conditions:
\begin{equation}
\label{psi_est}
(i)\ \ \Psi_{\an,\frac gf}\, (C)=
O\lr \Bigl( \frac{\log r_n}{r_n} \Bigr)^{1/4}\frac {|g|}f (\qa)\rr; \quad (ii)\ \ \widehat{\Psi}_{\an,\frac gf}\, (C)=
O\lr \Bigl( \frac{\log n}{r_n} \Bigr)^{1/4}\frac {|g|}f (\qa)\rr.
\end{equation}

We preface a~formulation  of the  corollaries of Theorems~\ref{stm_Bah_1}---\ref{stm_Bah_2} with a~stating  of two its direct consequence under  conditions~\eqref{psi_est}.
\begin{theorem}
\label{thm_2.3}
 Suppose that $r_n\to \infty$, as $\nty$, the condition $[A_1]$ holds true and $G$ is
differentiable on the set $F^{-1}(U)$. Assume in the addition  that the condition $(i)$ in ~\eqref{psi_est} holds true. Then the representation \eqref{stm_(A.1)} and  the relation \eqref{delta_(1.1)}  are valid together with $\Delta_n=A\, (\an(1-\an))^{1/4}\lr\frac{\log  r_n}n\rr^{3/4}\frac {|g|}f(\qa)$, where $A$ is some positive constant not depending on $n$.

Moreover, if additionally the condition $[A_2]$  and  relation $(ii)$ in ~\eqref{psi_est} are satisfied, then
\eqref{stm_(A.1)} and  \eqref{delta_h_(1.1)} are valid with $\widehat{\Delta}_n=A\, (\an(1-\an))^{1/4}\lr\frac{\log  n}n\rr^{3/4}\frac {|g|}f(\qa)$.
\end{theorem}
\begin{theorem}
\label{thm_2.4}
 Suppose that $r_n\to \infty$, as $\nty$, the condition $[A_1]$ holds true and $G$ is
differentiable on the set $F^{-1}(U)$. Assume in the addition  that the condition $(i)$ in ~\eqref{psi_est} holds true. Then the representation \eqref{stm_(A.2)} and  the relation \eqref{delta_(1.2)} are valid together with $\Delta_n=A\, (\an(1-\an))^{3/4}\lr\frac{\log r_n}n\rr^{5/4}\frac {|g|}f(\qa)$, where $A$ is some positive constant not depending on $n$.

Moreover, if additionally the condition $[A_2]$ and the relation  $(ii)$ in ~\eqref{psi_est} are  satisfied, then
\eqref{stm_(A.2)} and  \eqref{delta_h_(1.2)} are valid with $\widehat{\Delta}_n=A\, (\an(1-\an))^{3/4}\lr\frac{\log n}n\rr^{5/4}\frac {|g|}f(\qa)$.
\end{theorem}

Now we expose certain sets of conditions sufficient for the relations  \eqref{psi_est}  and obtain some  corollaries of  Theorems~\ref{thm_2.3}---\ref{thm_2.4}.

Let $SRV^{+\infty}_{\rho}$ ($SRV^{-\infty}_{\rho}$) be a~class of regularly varying in $+\infty$ ($-\infty$)
functions: \ \ $g\in SRV^{+\infty}_{\rho}$ ($SRV^{-\infty}_{\rho}$) $\Leftrightarrow$
 $(i)$ $g(x)=\pm |x|^{\rho}\, L(x)$, for $|x|>x_0$, with some $x_0>0$ ($x_0<0$),
$\rho\in\mathbb{R}$, and $L(x)$ is a positive slowly varying function at $+\infty$ ($-\infty$);  $(ii)$ the following second order regularity condition on the tails is satisfied
\begin{equation}
\label{r1} \lv \frac{}{}g(x+\triangle x)-g(x)\rv=O\Bigl( |\,g(x)|\,\,
\Big|\frac{\triangle x}{x}\, \Big|^{1/2}\Bigr),
\end{equation}
when $\triangle x=o(|x|)$, as $x\to +\infty$ ($x\to -\infty$).

Note that \eqref{r1} holds true for $g$ if $\Big|\frac
{L(x+\triangle x)}{L(x)}-1\Big|=O\Bigl( \Big|\, \frac{\triangle
x}{x}\, \Big|^{1/2}\Bigr)$, as $x\to +\infty$ ($x\to -\infty$), where $L$ is the
corresponding slowly varying function, and it is satisfied (even with degree $1$ instead of $1/2$) if $L$
is continuously differentiable for sufficiently large $|x|$ and
$|L'(x)|=O\lr\frac {L(x)}{|x|}\rr$, as $x\to +\infty$ ($x\to -\infty$), which is
valid for instance when $L$ is some power of the logarithm.
\begin{corollary} \label{cor3} Suppose that $p_n\to 0$ ($p_n\to 1$), condition $[A_1]$ is satisfied,  $f\in SRV_{\rho}^{-\infty}$ ($f\in SRV_{\rho}^{+\infty}$), \ where $\rho=-(1+\gamma)$, \ $\gamma >0$, and  $g\in SRV_{\rho}^{-\infty}$ ($g\in SRV_{\rho}^{+\infty}$), \ where $\rho\in\mathbb{R}$.
Then the condition  $(i)$ in  \eqref{psi_est} is satisfied, and if additionally $[A_2]$ holds true, then the condition  $(ii)$ (cf.~\eqref{psi_est}) is also satisfied. Hence, both  assertions stated in Theorems~\ref{thm_2.3}---\ref{thm_2.4} are valid.
\end{corollary}
We relegate the proof of  the Corollary \ref{cor3} to the  section~\ref{proof3}.

Our final corollary concerns the case when the $df$ $F$ and the function $G$ are twice differentiable.

Let us define a~function $\emph{v}\,(u)=\frac gf\circ F^{-1}(u)$, $u\in (0,1)$.

\begin{corollary} \label{cor4} Suppose that $p_n\to 0$ ($p_n\to 1$), condition $[A_1]$ is satisfied, and assume
  that the functions $f$, $g$ are differentiable on the set $F^{-1}(U)$. In the addition suppose that
\begin{equation}
\label{sup_1}
\sup_{u\in U}\lv \frac {\emph{v}\,'(u)\, [u\wedge (1-u)]}{\emph{v}\,(u)}\rv <\infty,
\end{equation}
and that
\begin{equation}
\label{sup_2}
\limsup_{u\downarrow 0\ (u\uparrow 1)}\lv \frac {\emph{v}\,\bigl(u+[u\wedge (1-u)]o(1)\bigr)}{\emph{v}\,(u)}\rv <\infty,
\end{equation}
where $o(1)$ denotes any function tending to zero when $u\downarrow 0$  $(u\uparrow 1)$.

Then the condition  $(i)$ in  \eqref{psi_est} is satisfied, and if additionally $[A_2]$ holds true, then the condition  $(ii)$ (cf.~\eqref{psi_est}) is also satisfied. Hence, both  assertions stated in Theorems~\ref{thm_2.3}---\ref{thm_2.4} are valid.
\end{corollary}
{\it Proof.} The proof of the corollary~\ref{cor4} is straightforward. Take an~arbitrary $C>0$,  fix $t:$ $|t|<C$, and put  $\alpha(n)=\sqrt{\frac{\log r_n}{r_n}}$ when we prove relation $(i)$ of  \eqref{psi_est}, and $\alpha(n)=\sqrt{\frac{\log n}{r_n}}$  when we prove relation $(ii)$ of  \eqref{psi_est} (under additional condition $[A_2]$). In both cases we have $\alpha(n)\to 0$, as $\nty$. Consider $\lv \emph{v}\,\bigl(p_n+t[p_n\wedge (1-p_n)]\alpha(n)\bigr)-\emph{v}\,(p_n)\rv$. Since for all sufficiently large $n$ \ $p_n$ and
 $p_n+t[p_n\wedge (1-p_n)]\alpha(n)$ belong to the set $U$, the latter quantity is equal
 \begin{equation}\label{proof_cor_2.4}
 \begin{split}
\lv \emph{v}\,(p_n)\rv\, \lv \frac{\emph{v}\,'\bigl(p_n+\theta t[p_n\wedge (1-p_n)]\alpha(n)\bigr)}{\emph{v}\,\bigl(p_n+\theta t[p_n\wedge (1-p_n)]\alpha(n)\bigr)}\,t\, [p_n\wedge (1-p_n)]\alpha(n)\rv \\
\times\lv\frac{\emph{v}\,\bigl(p_n+\theta t[p_n\wedge (1-p_n)]\alpha(n)\bigr)}{\emph{v}\,(p_n)} \rv=O\Bigr( \bigl| \emph{v}\,(p_n)\bigr| \alpha(n) \Bigr),
\end{split}
\end{equation}
what yields \eqref{psi_est}. The corollary is proved. $\quad \square$

\medskip
The following examples show that the conditions \eqref{sup_1} and \eqref{sup_2} hold true in a~number of interesting cases.
\begin{example}
 \label{ex_1}(Gumbel) Consider a distribution $F(x)=\exp(-\exp(-x))$, $x\in \mathbb{R}$, and let $g(x)=x^k$, where $k\in \mathbb{Z}=\{0,\pm 1,\pm 2,\dots \}$. We take  $k$ integer only to avoid some problems of the existence for negative $x$. In this case we have $f(x)=\exp(-x)\exp(-\exp(-x))$, for the inverse function we have $F^{-1}(u)=-\log(-\log u)$, $u\in (0,1)$. In this case we obtain $f(F^{-1}(u))=-u\, \log u$, and
$\emph{v}\,(u)=\frac{[-\log(-\log u)]^k}{-u\, \log u}$. After simple computations we obtain
\begin{equation}
\label{gumb}
\frac {\emph{v}\,'(u)\, [u\wedge (1-u)]}{\emph{v}\,(u)}=-k\frac{u\wedge (1-u)}{-\log(-\log u)\ u\, \log u }+\frac{u\wedge (1-u)}{-u\, \log u}(1+\log u).
\end{equation}
If $u\to 0$, the first term  at the r.h.s in \eqref{gumb} tends to zero and the second term tends to \ $-1$. When $u\to 1$ we obtain that the first term is equivalent $-k\frac{1-u}{-\log(-\log u)\, \log (1+(u -1)} \sim k\frac 1{log(-\log u)}=o(1)$. The second term is equivalent $-\frac{1-u}{\log u}=\frac{u-1}{\log (1+(u-1))}=1*o(1)$. Thus, \eqref{sup_1} is satisfied in both cases $U=[0,\varepsilon]$ ($p_n\to 0$) and $U=[1-\varepsilon,1]$ ($p_n\to 1$). The check \eqref{sup_2} we write $ \frac {\emph{v}\,\bigl(u+[u\wedge (1-u)]o(1)\bigr)}{\emph{v}\,(u)}=
\lb\frac{\log(-\log (u+[u\wedge (1-u)]o(1)))}{\log(-\log u)}\rb^k \frac u{u+[u\wedge (1-u)]o(1)}\frac{\log u}{\log(u+[u\wedge (1-u)]o(1))}$, and arguing as before we obtain that the latter quantity is $1+o(1)$, as $u\to 0$
and  as \ $u\to 1$ as well.
\end{example}
\begin{example} \label{ex_2} Let $F(x)=\lr 1-\exp(-x^{\gamma})\rr\mathbb{I}(x\ge 0) $, $\gamma >0$, and let $g(x)=x^{\rho}$, $\rho\in \mathbb{R}$. Now we get $F^{-1}(u)=[-\log (1-u)]^{1/\gamma}$, $u\in (0,1)$, and $\emph{v}\,(u)=\frac{[-\log (1-u)]^{\rho/\gamma}}{\gamma [-\log (1-u)]^{(\gamma-1)/\gamma}(1-u)} =\frac 1{\gamma}[-\log (1-u)]^{(\rho+1)/\gamma -1}\frac 1{1-u}$. Then we obtain
\begin{equation}
\label{gam}
\frac {\emph{v}\,'(u)\, [u\wedge (1-u)]}{\emph{v}\,(u)}=\frac{\rho+1-\gamma}{\gamma}\ \frac{u\wedge (1-u)}{-(1-u)\log(1-u)} + \frac{u\wedge (1-u)}{1-u}.
\end{equation}
The first term on the r.h.s in \eqref{gam} tends to the constant $\frac{\rho+1-\gamma}{\gamma}$ when $u\to 0$ and it tends to zero when $u\to 1$, the second term tends to zero, as $u\to 0$ and it tends to $1$, as  $u\to 1$. Thus,
\eqref{sup_1} is satisfied in both cases as in previous example. The check \eqref{sup_2} we write $ \frac {\emph{v}\,\bigl(u+[u\wedge (1-u)]o(1)\bigr)}{\emph{v}\,(u)}=\lb \frac{\log(1-u-[u\wedge (1-u)]o(1))}{\log(1-u)}\rb^{\frac{\rho+1}{\gamma}-1}\frac{1-u}{1-u-[u\wedge (1-u)]o(1)}$. The simple computations show that both factors of the latter quantity tends to $1$, as $u\to 0$ and as $u\to 1$.
\end{example}
\begin{example} \label{ex_3} (Weibull) Let $F(x)=\exp(-x^{-\gamma})\mathbb{I}(x\ge 0) $, $\gamma >0$, and let $g(x)=x^{\rho}$, $\rho\in \mathbb{R}$. Here we get  $F^{-1}(u)=[-\log u]^{-1/\gamma}$, $u\in (0,1)$, $f(F^{-1}(u))=
\gamma (-\log u)^{(\gamma+1)/\gamma}\, u$, and $\emph{v}\,(u)=\frac{[-\log u]^{-\rho/\gamma}}{\gamma [-\log u]^{(\gamma+1)/\gamma}\,u} =\frac 1{\gamma}[-\log (1-u)]^{-(\rho+\gamma+1)/\gamma }\frac 1{u}$. Then we obtain
\begin{equation}
\label{weib}
\frac {\emph{v}\,'(u)\, [u\wedge (1-u)]}{\emph{v}\,(u)}=-\frac{\rho+\gamma+1}{\gamma}\ \frac{u\wedge (1-u)}{u \, \log u} -\frac{u\wedge (1-u)}{u}.
\end{equation}
If  $u\to 0$, the  first term on the r.h.s in \eqref{weib} tends to zero and the second one tends to \ $-1$, and when $u\to 1$, the first term tends to the constant $\frac{\rho+\gamma +1}{\gamma}$ and the second one tends to zero. Thus, \eqref{sup_1} is satisfied in both cases  $u\to 0$,  $u\to 1$. The check \eqref{sup_2} we write $ \frac {\emph{v}\,\bigl(u+[u\wedge (1-u)]o(1)\bigr)}{\emph{v}\,(u)}=\lb \frac{\log u}{\log(u+[u\wedge (1-u)]o(1))}\rb^{\frac{\rho+\gamma +1}{\gamma}}\frac{u}{u+[u\wedge (1-u)]o(1)}$, and simple evident arguments  show that both factors here  tends to $1$, as $u\to 0$ and as $u\to 1$.
\end{example}
\begin{example}  \label{ex_4}Let $C_{\gamma}\exp(-|x|^{\gamma})$, $\gamma>0$, where $C_{\gamma}$ is a~constant, depending only on
$\gamma$, and let $g(x)=\pm |x|^{\rho}$, $\rho\in \mathbb{R}$. It is clear that the asymptotic behavior of the functions at the l.h.s.'s in conditions \eqref{sup_1} and \eqref{sup_2} are similar as in example~\ref{ex_2} ($u\to 1$). So, these conditions are also satisfied.
\end{example}
\begin{example}  \label{ex_5} Here we consider an~example of a~distribution with super heavy tails, having  no finite moments. In this case some difficulties arise, nevertheless  the  Bahadur -- Kiefer representations \eqref{stm_(A.1)}--\eqref{stm_(A.2)} are  still valid for the intermediate sample  quantiles under some additional conditions.

Let $F(x)=1-\frac C{\log x}$ for $x\ge x_0>0$, where $C>0$ is some constant. Suppose for ease of presentation that $p_n\to 1$, as $\nty$, while $r_n=n-k_n\to \infty$, and let $g(x)=x^{\rho}$, $\rho\in \mathbb{R}$, though this will not influence the basic outline of our results.

In this case $F^{-1}(u)=\exp\lr\frac C{1-u}\rr$, $f(F^{-1}(u))=\frac{(1-u)^2}C \exp\lr-\frac C{1-u}\rr$,
$\emph{v}\,(u)=\exp\lr (\rho +1) \frac C{1-u}\rr\frac C{(1-u)^2}$. Since $p_n\to 1$, we are interested only in the case $u\to 1$, so \,$u\wedge (1-u)=1-u$, and after simple computations we obtain
\begin{equation}
\label{heavy}
\frac {\emph{v}\,'(u)\, (1-u)}{\emph{v}\,(u)}=\frac{C(\rho+1)}{1-u} +2,
\end{equation}
what is not bounded as $u\to 1$, and therefore~\eqref{sup_1} is clearly not satisfied. The computations of the magnitude on the l.h.s. in  \eqref{sup_2} yields
\begin{equation}
\label{heavy1}
\frac {\emph{v}\,\bigl(u+[u\wedge (1-u)]o(1)\bigr)}{\emph{v}\,(u)}=\exp\lr C(\rho +1)\frac{o(1)}{1-u}\rr \lr 1+o(1)\rr.
\end{equation}
We conclude that  \eqref{sup_2} is satisfied only if $\frac{o(1)}{1-u}\to 0$, as $u\to 1$. However,  we apply our conditions
for a~sequence with  $u=p_n$ (cf.~proof of the Corollary~\ref{cor4}). So, $1-u=1-p_n$, the quantity $o(1)$ is $\alpha(n)=\sqrt{\frac {\log r_n}{r_n}}$ (cf.~\eqref{proof_cor_2.4}), where $r_n=k_n \wedge (n-k_n)$ (and $r_n=n-k_n$ for all sufficiently large $n$,as $p_n\to 1$).
 Although the relations  \eqref{psi_est} are not valid more in our example, we can achieve a~weaker relation $ \Psi_{\an,\frac gf}\, (C)=o\bigl(\frac {|g|}f(\qa)\bigr)$, guaranteeing that \eqref{stm_(A.1)} and \eqref{stm_(A.2)} are representations (cf.~Remark~\ref{rem_1}). Observe that \eqref{proof_cor_2.4} (cf.~the proof of Corollary~\ref{cor4}) and  \eqref{heavy}--\eqref{heavy1} together imply that we need only that  $\frac{\alpha(n)}{1-p_n}=
 \sqrt{\frac {\log r_n}{r_n}}\ \frac n{r_n}=o(1)$, as $\nty$. Thus, the representations~\eqref{stm_(A.1)} and \eqref{stm_(A.2)} are valid for the intermediate sample quantiles in this  example if \
 \begin{equation*}
 \frac{n^{2/3}(\log r_n)^{1/3}}{r_n}\to 0, \text{ \ \ as } \quad \nty.
 \end{equation*}
\end{example}

Define a~binomial r.v.  $\na=\sharp \{i \ : X_i \le \xia \}$, $0<\alp <1$.  Our proof of Theorems~\ref{stm_Bah_1}-\ref{stm_Bah_2}  uses the following fact: conditionally on $\na$ the
order statistics $X_{1:n},\dots ,X_{\na:n}$ are distributed as
order statistics corresponding to a~sample of $\na$ i.i.d.~r.v.'s
with distribution function $F(x)/\alp$, $x\le \xia$. Though this
fact is well known (cf., e.g.,~Theorem 12.4,~\cite{kallen}, cf.~also
\cite{gh2007}, \cite{hp}), we give a brief proof of it in the section~\ref{appendix}.
 \section{Proof of Theorem~\ref{stm_Bah_1}}
\label{proof1}
We can assume with impunity that $a_2\le 1/2$, i.e. we will prove representation~\eqref{stm_(A.1)} for
the quantiles at the left edge of the variation series. Then $\kn \le (n-\kn)$ for all sufficiently large $n$, and so it is enough to prove \eqref{stm_(A.1)} with
\begin{equation}
\label{delta1}
\Delta_n=A\, \an^{1/4}\lr\frac{\log k_n}n\rr^{3/4}\frac {|g|}f(\qa)+ B\,\an^{1/2}\lr\frac{\log k_n}n\rr^{1/2}  \Psi_{\an,\frac gf}\, (C).
\end{equation}
We begin with the proof of the first assertion of the theorem, where there is no restrictions on $k_n$ in its tending to infinity.

Let $U_1,\dots,U_n$ denote a sample of independent uniform $(0,1)$
distributed r.v.'s, and $U_{1:n}\le\cdots\le U_{n:n}$ -- the
corresponding order statistics. Put
\begin{equation}
\label{stm_(A.3)}
N_{\an}^x = \sharp \{i : X_i \le \xi_{\an} \}\, ,\quad N_{\an} = \sharp
\{i : U_i \le \an \},
\end{equation}
and note that $\qan =X_{\kn:n}$ (because $\an=\kn/n$).

We must prove that
$\boldsymbol{P}(|R_n(p_n)|>\Delta_n)=O\lr k_n^{-c}\rr$ for each $c>0$ (cf. \eqref{stm_(A.1)}),
and since the joint distribution of $X_{\kn :n}$, $N_{\an}^x $
coincide with joint distribution of $F^{-1}(U_{\kn :n})$,
$N_{\an}$ it is suffices to verify it for a remainder given by
\begin{equation*}
R_n(p_n)= G(F^{-1}(U_{\kn :n}))-G(F^{-1}(\an))+\frac{N_{\an}-\an
n}{n}\frac gf(\qa).
\end{equation*}
Since $\boldsymbol{P}(U_{\kn :n}\notin U)=O(exp(-\delta n))$ for some
$\delta >0$ not depending on $n$, we can rewrite $R_n(p_n)$ for all
sufficiently large $n$ as
\begin{equation}
\label{stm_(A.4)}
\frac gf(\qa)\, R_{n,1}+ R_{n,2},
\end{equation}
where $R_{n,1}=U_{\kn :n}-\an+\frac{N_{\an}-\an n}{n}$, and
$R_{n,2}=\Bigl( \frac gf\lr F^{-1}(\an+\theta(U_{\kn :n}-\an
))\rr$ $-\frac gf \lr F^{-1}(\an)\rr\Bigr)\bigl(U_{\kn
:n}-\an\bigr)$, $0<\theta <1$. Fix an arbitrary $c>0$ and note
that we can estimate $R_{n,j}$, $j=1,2$, on the set $E=\bigl\{
\omega : |N_{\an}-\an n|<A_0\bigl(\an \, n\, \log k_n \bigr)^{1/2}
\bigr\}$, where $A_0$ is a positive constant, depending only on
$c$, because by Bernstein inequality $\boldsymbol{P}(\Omega \setminus
E)=O(k_n^{-c})$ (in fact we can take every $A_0$:  $A^2_0>2c$). We will prove that
\begin{equation}
 \label{stm_(A.5)}
\boldsymbol{P}\bigl( |R_{n,1}|>A_1(\, \an\, )^{1/4}(\log k_n/n)^{3/4}\bigr)
=O(k_n^{-c})
\end{equation}
and that
\begin{equation}
\label{stm_(A.6)}
\boldsymbol{P}\bigl( |R_{n,2}|>A_2\, \an\, \Psi_{\an,\frac
gf}(C)(\log k_n/\kn)^{1/2}\bigr) =O(k_n^{-c}).
\end{equation}
Here and elsewhere $A_i\, , i=1,2,\dots $, and $C$ denote some positive
constants, depending only on $c$. Relations \eqref{stm_(A.4)}--\eqref{stm_(A.6)} imply \eqref{stm_(A.1)} with $\Delta_n$ given in \eqref{delta1}.

First we prove \eqref{stm_(A.5)}, using a similar conditioning on $N_{\an}$
argument as in proof of lemmas 4.1, 4.3 in~\cite{gh2007}. First let $\kn\le N_{\an}$, then conditionally on
$N_{\an}$ the order statistic $U_{\kn:n}$ is distributed as
$\kn$-th order statistic $U'_{\kn:N_{\an}}$ of the sample
$U'_1,\dots,U'_{N_{\an}}$ independent $(0,\an)$ uniformly
distributed r.v.'s (cf. lemma~\ref{Lemma_A}, appendix). Its expectation
$\boldsymbol{E}\big(U_{\kn:n}\lv\rp N_{\an},\ \kn\le
N_{\an}\big)=\an\, \frac{\kn}{N_{\an}+1}$, and the conditional
variance $V^2_{\kn}=\frac{\alp_n^2}{N_{\an}+2}\,
\frac{\kn}{N_{\an}+1}\, \bigl( 1-\frac{\kn}{N_{\an}+1}\bigr)$, and
on the set $E$ we have an estimate $V^2_{\kn}\le A_0(\an)^{1/2}\,
n^{-3/2}\log^{1/2}n$. \ \ Then rewrite $R_{n,1}$ (at the event $\kn\le N_{\an}$) as
\begin{equation}
\label{stm_(A.7)}
U_{\kn:n}-\an\, \frac{\kn}{N_{\an}+1}+R'_{n,1},
\end{equation}
where $R'_{n,1}=\an\, \frac{\kn}{N_{\an}+1}-\an+\frac{N_{\an}-\an
n}{n}=\frac{(N_{\an}-\kn)^2}{n(N_{\an}+1)}+\frac{N_{\an}-\kn}{n(N_{\an}+1)}
-\frac{\kn}{n(N_{\an}+1)}$, and on the set $E$ the latter
quantity is of the order $O\lr\frac{\log k_n}{n}\rr$, and since
$\frac{\log k_n}{n}=o\bigl( \an^{1/4}\lr
\frac{\log k_n}{n}\rr^{3/4}\bigr)$, the remainder term $R'_{n,1}$
is of negligible order for our purposes. For the first two terms
in \eqref{stm_(A.7)} we have
\begin{equation*}
\boldsymbol{P}\lr \lv U_{\kn:n}-\an\, \frac{\kn}{N_{\an}+1}\rv>A_1
(\an)^{1/4}\lr \frac{\log k_n}{n}\rr^{3/4}\big| N_{\an}:\ \kn\le
N_{\an}\rr
\end{equation*}
\begin{equation}
\label{stm_(A.8)}
=\boldsymbol{P}\lr \lv U'_{\kn:N_{\an}}-\an\, \frac{\kn}{N_{\an}+1}\rv>A_1
(\an)^{1/4}\lr \frac{\log k_n}{n}\rr^{3/4}\rr =P_1+P_2, \quad
\end{equation}
where $N_{\an}$ is fixed, $\kn\le
N_{\an}$, $A_1$ is a constant which we will choose later,
$P_1=\boldsymbol{P}\lr U'_{\kn:N_{\an}}> \an\, \frac{\kn}{N_{\an}+1}+A_1
(\an)^{1/4}\lr \frac{\log k_n}{n}\rr^{3/4}\rr$,
$P_2=\boldsymbol{P}\lr U'_{\kn:N_{\an}}< \an\, \frac{\kn}{N_{\an}+1}-A_1
(\an)^{1/4}\lr \frac{\log k_n}{n}\rr^{3/4}\rr$.
We evaluate $P_1$, the treatment for $P_2$ is similar. Consider a
binomial r.v.
$S'_n=\sum_{i=1}^{N_{\an}}\bf{1}_{\{U'_{i:N_{\an}}\le \an\,
\frac{\kn}{N_{\an}+1}+A_1 (\an)^{1/4}\lr
\frac{\log k_n}{n}\rr^{3/4} \}}$ with parameter $(q_n,N_{\an})$,
where $q_n=\min\bigl(1,\frac{\kn}{N_{\an}+1}+t_n \bigr)$, where $t_n=A_1 \lr
\frac{\log k_n}{\kn}\rr^{3/4} $. If $q_n=1$, then $P_1=0$
and the inequality we need is valid trivial. Let $q_n<1$ and let $\overline{S'}_n$ denote the average $S'_n/N_{\an}$, then  the probability  $P_1$ is equal to
\begin{equation}
\label{stm_(A.9)}
\boldsymbol{P}(S'_n<\kn)= \boldsymbol{P}\lr\overline{S'}_n-q_n<\frac{k_n}{N_{\an}}-\frac{k_n}{N_{\an}+1} -t_n
\rr.
\end{equation}
Note that $\frac{k_n}{N_{\an}}-\frac{k_n}{N_{\an}+1}=\frac{k_n}{N_{\an}(N_{\an}+1)}<\frac 1{N_{\an}}$, and since  the latter quantity is $o\lr t_nk_n^{-1/4}\rr=o(t_n)$ on the set $E$,  this term  can be omitted at the r.h.s. of
\eqref{stm_(A.9)} in our estimating. To evaluate  $ \boldsymbol{P}\lr\overline{S'}_n-q_n< -t_n\rr$
we note that $q_n-t_n=\frac{k_n}{N_{\an}+1}\in (0,1)$, and that $q_n>1/2$ for all sufficiently large $n$ (and hence $k_n$ and  $N_{\an}$) on the set $E$. So, we may apply an inequality (2.2) of Hoeffding \cite{hoef} with $\mu=q_n$ and with $g(\mu)=1/(2\mu(1-\mu))$. Then we obtain
\begin{equation}
\label{stm_(A.9_h)}
\boldsymbol{P}(S'_n<\kn)\le \exp \lr-N_{\an}t_n^2g(q_n)\rr = \exp \lr-\frac{N_{\an}A_1^2 \bigl(\log k_n/\kn\bigr)^{3/2}}{2q_n(1-q_n)}\rr.
\end{equation}
Finally we note that  $1-q_n=1-\frac{\kn}{N_{\an}+1}-A_1 \lr
\frac{\log k_n}{\kn}\rr^{3/4}\le \frac{N_{\an}+1-\kn}{N_{\an}+1}$, and on the set $E$ the latter quantity is not greater than $\frac{A_0(\kn\log k_n)^{1/2}}{N_{\an}}$. Then we can get a low bound for the ratio at the r.h.s. in \eqref{stm_(A.9_h)}: $\frac{N_{\an}A_1^2 \bigl(\log k_n/\kn\bigr)^{3/2}}{2q_n(1-q_n)} \ge \frac{A_1^2 N^2_{\an} \bigl(\log k_n/\kn\bigr)^{3/2}}{2A_0(\kn\log k_n)^{1/2}}=\frac{A_1^2}{2A_0} \log k_n\lr\frac{N_{\an}}{k_n}\rr^2=\frac{A_1^2}{2A_0} \log k_n\lr 1+o(1)\rr$. This bound and \eqref{stm_(A.9_h)} together
yield that when $\frac{A_1^2}{2A_0}\ge c$ the desired  relation  $P_1=O(k_n^{-c})$ hold true. The same estimate is valid for $P_2$.

Note that the condition
$\frac{A_1^2}{2A_0}\ge c$ which we  needed to establish the desired estimates can be weakened to $\frac{A_1^2}{2A_0}\ge c-1/2>0$ if we apply a refinement of Heoffding's inequality due to Talagrand \cite{Talagrand} (cf.~also Leon and Perron \cite{Leon_03}). However the improvement is not very useful here, as applying Talagrand's inequality instead of Hoeffding,s only affects the constant, but not the order bound in our setting.

In case $N_{\an}<\kn$ we use the fact that $U_{\kn:n}$
conditionally on $N_{\an}$ is distributed as $(\kn-N_{\an})$-th
order statistic $U''_{\kn-N_{\an}:n-N_{\an}}$ of the sample
$U''_1,\dots,U''_{n-N_{\an}}$ from $(1-\an,1)$ uniform
distribution, its expectation is
$\an+\frac{\kn-N_{\an}}{n-N_{\an}+1}$, and for the conditional
variance we have the estimate $V^2_{\kn-N_{\an}}\le
A_0( \an\log k_n\,)^{1/2}\, n^{-3/2}$. In this case we use a
representation for $R_{n,1}=R''_{n,1}+R''_{n,2}$, where
$R''_{n,1}=U_{\kn:n}-\an-\frac{\kn-N_{\an}}{n-N_{\an}+1}(1-\an)$,
and $R''_{n,2}=\frac{N_{\an}-\an\,
n}{n}+\frac{\kn-N_{\an}}{n-N_{\an}+1}(1-\an)$. Similarly as in
first case we obtain that
$R''_{n,2}=O\bigl(\frac{\log k_n}{n}\bigr)$ with probability
$1-O(k_n^{-c})$, and this term is of the negligible order in our
estimating. Using Hoeffding's  inequality we obtain for $R''_{n,1}$
same estimate as for $R'_{n,1}$. So \eqref{stm_(A.5)} is proved.

It remains to prove  \eqref{stm_(A.6)}. First note that by  \eqref{stm_(A.5)}
on the set $E$
with probability $1-O(k_n^{-c})$ we have $|U_{\kn:n}-\an|\le
A_0\frac{(\kn\log k_n)^{1/2}}{n}+A_1\an\lr\frac{\log k_n}{\kn}\rr^{3/4}
=\bigl(\an\frac{\log k_n}{n}\bigr)^{1/2}\bigl( 1+o(1)\bigr)$. Thus,
there exists $A_2$, depending only on $c$, such that $|R_{n,2}|\le
A_2 \bigl(\an\frac{\log k_n}{n}\bigr)^{1/2}\Psi_{\an,\frac
gf}(A_2)$ with probability $1-O(k_n^{-c})$. This implies \eqref{stm_(A.6)}.
Thus, the first assertion of the theorem~\ref{stm_Bah_1} is proved.

To prove the second assertion, it is sufficient to repeat previous arguments replacing $\log k_n$ by $\log n$ throughout the proof, and applying the fact $\log n/k_n\to 0$ (due to $[A_2]$)  instead of the evident fact that $\log k_n/k_n\to 0$ used before, moreover now we should use the function $\widehat{\Psi}_{\an,h}(C)$ instead of $\Psi_{\an,h}(C)$. These replacements lead to estimates with probability $O(n^{-c})$ for each $c>0$. The theorem is proved. $\quad \square$
\section{Proof of theorem~\ref{stm_Bah_2}}
\label{proof2}
We give a~detailed proof of the first assertion of Theorem~\ref{stm_Bah_2}. To prove the second one it is enough to  make similar replacements as in the proof of the corresponding part of the Theorem~\ref{stm_Bah_1}, therefore we omit it.

Let $N^x_{\an}$ and $N_{\an}$
are given as in \eqref{stm_(A.3)}, then we can rewrite integral on the l.h.s. of  \eqref{stm_(A.2)}
as $\frac {sgn(N^x_{\an}-\kn )}{n}\sum_{i=(\kn\wedge N^x_{\an})+1}^{\kn\vee N^x_{\an}}(G(X_{i:n})-G(\qa))$, where $sgn(x)=x/|x|$, $sgn(0)=0$. Let us adopt the following notation: for
any integer $k$ and $m$ define a~set $I_{(k,m)}:=\{i: (k\wedge m)+1\le i\le k\vee m\}$ and let
$\sum_{i\in I_{(k,m)}}(.)_i:=sgn(m-k)\sum_{i=(k\wedge m)+1}^{k\vee m}(.)_i$. Then we must estimate
$R_n(p_n)=\frac
1n\sum_{i\in I_{(\kn,{N^x_{\an})}}}(G(X_{i:n})-G(\qa))+\frac{\bigl(N^x_{\an}-\an n\bigr)^2}{2n^2}\frac gf(\qa)$ \ (cf.~ \eqref{stm_(A.2)}), and
similarly as in proof of Theorem~\ref{stm_Bah_1} we note that $R_n(p_n)$ is
distributed as
\begin{equation*}
\frac 1n\sum_{i\in I_{(\kn,{N_{\an})}}}\Bigl(G\circ F^{-1}(U_{i:n})-G\circ
F^{-1}(\an)\Bigr)+\frac{\bigl(N_{\an}-\an n\bigr)^2}{2n^2}\frac
gf(\qa) \phantom{R_{n,1}+R_{n,2}}
\end{equation*}
\begin{equation}
\label{stm_(A.10)}
\phantom{R_{n,1}+R_{n,2}R_{n,1}+R_{n,2}+R_{n,1}}=\frac gf(\qa)R_{n,1}+R_{n,2},
\end{equation}
where \\ $R_{n,1}=\frac 1n\sum_{i\in I_{(\kn,{N_{\an})}}}(U_{i:n}-\an)+\frac{\bigl(N_{\an}-\an n\bigr)^2}{2n^2}$\, ,\\
$R_{n,2}=\frac 1n\sum_{i\in I_{(\kn,{N_{\an})}}}\lb \frac gf\circ
F^{-1}\bigl(\an +\theta_i(U_{i:n}-\an)\bigr)-\frac gf\circ
F^{-1}\bigl(\an\bigr)\rb\Bigl(U_{i:n}-\an\Bigr)$, \\
è \ $0<\theta_i<1$, \ \ $i\in I_{(\kn,{N_{\an})}}$.

As well as before (cf. the proof of Theorem~\ref{stm_Bah_1}) we can assume with impunity that $a_2\le 1/2$, then we need to prove \eqref{stm_(A.2)} with
\begin{equation}
\label{delta2}
\Delta_n=A\, \an^{3/4}\lr\frac{\log k_n}n\rr^{5/4}\frac {|g|}f(\qa)+ B\,\an\frac{\log k_n}n  \Psi_{\an,\frac gf}\, (C),
\end{equation}

Fix an arbitrary $c>0$ and prove that
\begin{equation}
\label{stm_(A.11)}
\boldsymbol{P}\Bigl( |R_{n,1}|>A_1(\, \an\, )^{3/4}(\log k_n/n)^{5/4}\Bigr)
=O(k_n^{-c}),
\end{equation}
\begin{equation}
\label{stm_(A.12)}
\boldsymbol{P}\Bigl( |R_{n,2}|>A_2\, \an\, \frac{\log k_n}{n}\Psi_{\an,\frac
gf}(A_2)\Bigr) =O(k_n^{-c}),
\end{equation}
where  $A_i>0$, $i=1,2,\dots$,  are some  constants, depending only on $c$. Relations \eqref{stm_(A.10)} and \eqref{stm_(A.11)}--\eqref{stm_(A.12)} imply
\eqref{stm_(A.2)} with $\Delta_n$ as in \eqref{delta2}.
Similarly as when proving of Theorem~\ref{stm_Bah_1} it is enough to estimate
$R_{n,j}$, $j=1,2$, on the set $E=\bigl\{ \omega : |N_{\an}-\an
n|<A_0\bigl(\an \, n\, \log k_n \bigr)^{1/2} \bigr\}$, where
$A_0>0$ is a constant, depending only on $c$, such that $\boldsymbol{P}(\Omega
\setminus E)=O(k_n^{-c})$.

First we treat $R_{n,2}$. Note that
\begin{equation*}
\max_{i\in I_{(\kn,{N_{\an})}}}\big|\, U_{i:n}-\an\big|=\big|\,
U_{\kn:n}-\an\big|\vee \big|\, U_{N_{\an}:n}-\an\big|\vee \big|\, U_{N_{\an}+1:n}-\an\big|\, ,
\end{equation*}
$\boldsymbol{P}\Bigl(\big|\, U_{\kn:n}-\an\big|>A_0\bigl(\an\log k_n /n\,
\bigr)^{1/2}\Bigr)=O(k_n^{-c})$ (cf. proof of Theorem~\ref{stm_Bah_1}), and for
$j=N_{\an:n}\, , N_{\an:n}+1$ simultaneously we have  $\boldsymbol{P}\Bigl(\big|\,
U_{j:n}-\an\big|>A_1\frac{\log k_n}{n}\Bigr)\le \boldsymbol{P}\Bigl(
U_{N_{\an}+1:n}-U_{N_{\an}:n}>A_1\frac{\log k_n}{n}\Bigr) =\boldsymbol{P}\Bigl(U_{1:n}>A_1\frac{\log k_n}{n}\Bigr)=
\Bigl(1-A_1\frac{\log k_n}{n}\Bigr)^n=O(k_n^{-c})$ for $A_1>c$.
Since $\frac{\log k_n}{n}=o(\frac{\an\log k_n}{n})^{1/2}$, \ on the
set $E$ we obtain
\begin{equation*}
\big|R_{n,2}\big|\le \frac 1n\Psi_{\an\frac
gf}(A_0)A_0^2\Bigl(\an\, n\log k_n\Bigr)^{1/2}\Bigl(\frac{\an\,
\log k_n}{ n}\Bigr)^{1/2}=A_2\an\frac{\log k_n}{n}\Psi_{\an\frac
gf}(A_0)
\end{equation*}
with probability  $1-O(k_n^{-c})$, and   \eqref{stm_(A.12)} is proved.

Finally, consider $R_{n,1}$. Note that conditionally on $N_{\an}$,
\ $\kn < N_{\an}$, the order statistics $U_{i:n}$, $\kn\le i\le N_{\an}$, are distributed
as the order statistics $U'_{i:N_{\an}}$ from the uniform $(0,\an)$ distribution
(cf.~proof of theorem~\ref{stm_Bah_1}),  their conditional
expectations are equal to $\an\frac{i}{N_{\an}+1}$. Then in the case $\kn <  N_{\an}$ (the proof for the~case $N_{\an}\ge \kn$ is similar
(cf. proof of theorem~\ref{stm_Bah_1}) with respect to interval $(1-\an,1)$, and
we omit the details) we rewrite $R_{n,1}$ as
\begin{equation}
\label{stm_(A.13)}
R_{n,1}=\frac
1n\sum_{i=\kn+1}^{N_{\an}}\bigl(U_{i:n}-\an\frac{i}{N_{\an}+1}\bigr)+R'_{n,1}\,
,
\end{equation}
where $R'_{n,1}=\frac
1n\sum_{i=\kn+1}^{N_{\an}}\an\bigl(\frac{i}{N_{\an}+1}-1\bigr)+ \frac{(N_{\an}-\an\,
n)^2}{2n^2}=-\frac{\kn}{n^2}\frac{(N_{\an}-\kn)(N_{\an}-\kn-1)}{2\,(N_{\an}+1)}+ \frac{(N_{\an}-\kn)^2}{2n^2}
=\frac{(N_{\an}-\kn)^2(N_{\an}+1-\kn)}{2\,(N_{\an}+1)n^2}-\frac{\kn(N_{\an}-\kn)}{2(N_{\an}+1)n^2}$, and on the set $E$ the latter quantity is of the order $O\lr
\frac{\kn^{1/2}(\log k_n)^{3/2}}{n^2}\rr=o\lr\bigl(\an
\bigr)^{3/4}\lr\frac{\log k_n}{n}\rr^{5/4}\rr$,  i.e. $R'_{n,1}$ is
of negligible order (cf.~\eqref{stm_(A.11)}) for our purposes.

It remains to evaluate the dominant first term on the r.h.s. in \eqref{stm_(A.13)}.
Fix an arbitrary $c_1>c+1/2$, and note that conditional on
$N_{\an}$ the variance of $U_{i:n}$ ($\kn+1\le i\le N_{\an}$) is equal
to $V_i^2=\bigl(\an\bigr)^2\frac 1{N_{\an}+2}\,
\frac{i}{N_{\an}+1}\, \lr 1-\frac{i}{N_{\an}+1}\rr$, and on the
set $E$ it is less than
$\bigl(\an\bigr)^2\frac{A_0\kn^{1/2}(\log k_n)^{1/2}}{N^2_{\an}}$,
and $V_i\le\an A_0^{1/2}\kn^{1/4}(\log k_n)^{1/4}/N_{\an}
\le A_0^{1/2}\an \kn^{-3/4}(\log k_n)^{1/4}\le A_0^{1/2}\bigl(\an\bigr)^{1/4}n^{-3/4}(\log k_n)^{1/4}$. Using
Hoeffding's  inequality (similarly as in proof of theorem~\ref{stm_Bah_1}), we find
that
\begin{equation*}
\boldsymbol{P}\Bigl( \big|U_{i:n}-\an\frac{i}{N_{\an}+1} \big|>A_1
\bigl(\an\bigr)^{1/4}\bigl(\log k_n/n \bigr)^{3/4}\Big| N_{\an}:\
\kn\le N_{\an}\Bigr)=O(k_n^{-c})\, ,
\end{equation*}
where $A_1$ depends only on $c_1$ (in fact it is true for every
$A_1$ such that $A^2_1>2A_0c$). Thus
\begin{equation*}
\boldsymbol{P}\Bigl( \frac 1n
\big|\sum_{i=\kn}^{N_{\an}}\bigl(U_{i:n}-\an\frac{i}{N_{\an}+1}
\bigr)\big|>A_0A_1 \bigl(\an\bigr)^{3/4}\bigl(\log k_n/n
\bigr)^{5/4}\Big| N_{\an}:\ \kn\le N_{\an}\Bigr)
\end{equation*}
\begin{equation}
\label{stm_(A.14)}
\qquad \qquad \qquad \qquad \qquad
\qquad\qquad \quad \ \le A_0(\kn\log k_n )^{1/2}\,
O(k_n^{-c_1})=O(k_n^{-c})\, .
\end{equation}
Combining \eqref{stm_(A.13)}--\eqref{stm_(A.14)} and similar estimates for the case
$N_{\an}<\kn$, arrive at \eqref{stm_(A.11)}. The theorem is proved.$\quad \square$

\section{Proof of Corollary~\ref{cor3}}
\label{proof3}
Suppose for definiteness that $p_n\to 0$, as $\nty$, then we must prove the relations:
\begin{equation}
\label{psi_est1}
\Psi_{\an,\frac gf}\, (C)=
O\lr \Bigl( \frac{\log k_n}{k_n} \Bigr)^{1/4}\frac {|g|}f (\qa)\rr; \quad \widehat{\Psi}_{\an,\frac gf}\, (C)=
O\lr \Bigl( \frac{\log n}{k_n} \Bigr)^{1/4}\frac {|g|}f (\qa)\rr.
\end{equation}
Let $\log(\cdot)$ denote
$\log k_n$ when we prove a first of relations \eqref{psi_est1} and  $\log n$ when we prove  the second one.
Since we will need only that $\log(\cdot)/k_n\to 0$, what is evident in the first case and is valid  by
$[A_2]$ in the second one, this notation will allow us to prove each of desired assertions simultaneously.

Define $x_n=F^{-1}(p_n)$, which tend to $-\infty$, as $\nty$. Fix $C>0$ and for a~fixed $t:\ |t|\le C$, \ put \ $\triangle x_n=F^{-1}\lr p_n+t\sqrt{p_n \frac{\log(\cdot)}{n}}\rr-x_n=F^{-1}\lr p_n \Bigl(1+t\sqrt{ \frac{\log(\cdot)}{k_n}}\Bigr)\rr-x_n$.

First we prove that $\frac{\triangle x_n}{x_n}\to 0$, as $\nty$. Due to smoothness condition $[A_1]$ for all sufficiently large $n$ we may write $\frac{\triangle x_n}{x_n}=\frac{1}{x_n\, f\Bigl(F^{-1}\lr p_n \bigl(1+\theta \, t\sqrt{ \frac{\log(\cdot)}{k_n}}\bigr)\rr\Bigr)}t\,\sqrt{p_n \frac{\log(\cdot)}{n}}
=\frac 1{x_n f(x_n)}\frac{f(F^{-1}(p_n))}{f\Bigl(F^{-1}\lr p_n \bigl(1+\theta \, t\sqrt{ \frac{\log(\cdot)}{k_n}}\bigr)\rr\Bigr)}
t\,\sqrt{p_n \frac{\log(\cdot)}{n}}$, where $0<\theta<1$, and since due to regularity condition we have $f(x_n)x_n \sim -\gamma \, F(x_n)=-\gamma \, p_n$, as $x_n\to -\infty$ (cf., e.g.,~Bingham et al.~\cite{bingham}), the latter quantity is equivalent to $- \frac 1{\gamma \, p_n}\frac{f(F^{-1}(p_n))}{f\Bigl(F^{-1}\lr p_n \bigl(1+o(1)\bigr)\rr\Bigr)}
t\,\sqrt{p_n \frac{\log(\cdot)}{n}}=- \frac 1{\gamma }\frac{f(F^{-1}(p_n))}{f\Bigl(F^{-1}\lr p_n \bigl(1+o(1)\bigr)\rr\Bigr)}
t\,\sqrt{\frac{\log(\cdot)}{k_n}}$. It remains to show that $\frac{f(F^{-1}(p_n))}{f\Bigl(F^{-1}\lr p_n \bigl(1+o(1)\bigr)\rr\Bigr)}=1+o(1)$.
Since $f\in SRV^{-\infty}_{-(1+\gamma)}$, for all $x<x_0<0$ we have $f(x)=|x|^{-(1+\gamma)}L(x)$, where $L(x)$ is a~slowly varying in $-\infty$ positive function. Moreover, the inverse function $F^{-1}(u)$ is regular varying at zero, i.e. $F^{-1}(u)=u^{-1/\gamma}L_1(u)$, where $L_1(u)$ is a correspondent slowly varying at zero function. So, for sufficiently large $n$ we have $\frac{f(F^{-1}(p_n))}{f\Bigl(p_n \bigl(1+o(1)\bigr)\Bigr)}=\frac{\lb p_n^{-1/\gamma}L_1(p_n)\rb^{-(1+\gamma)}L(F^{-1}(p_n))}{\lb \bigl(p_n(1+o(1))\bigr)^{-1/\gamma}L_1(p_n(1+o(1)))\rb^{-(1+\gamma)}L(F^{-1}(p_n(1+o(1))))}\sim \frac{L(F^{-1}(p_n))}{L(F^{-1}(p_n(1+o(1))))}=\frac{L\lb p_n^{-1/\gamma}L_1(p_n)\rb}{L\lb\bigl(p_n(1+o(1))\bigr)^{-1/\gamma}L_1(p_n(1+o(1)))\rb}\sim 1$. Thus, $\lv\frac{\triangle x_n}{x_n}\rv =O\lr\sqrt{\frac{\log(\cdot)}{k_n}}\rr$.

Finally, we obtain a~bound for $\lv\frac gf\bigl(x_n+\triangle x_n\bigr)-\frac gf\bigl(x_n\bigr)\rv$ for an~arbitrary  fixed $C>0$ and $|t|\le C$, as $\nty$. Due to relation \eqref{r1} which holds true for the density $f$ as well as for the function $g$ we have
\begin{equation*}
\label{sec_4_2}
\begin{split}
&\lv\frac gf\bigl(x_n+\triangle x_n\bigr)-\frac gf\bigl(x_n\bigr)\rv=\lv\frac{f(x_n)[g(x_n+\triangle x_n)-g(x_n)]-g(x_n)[f(x_n+\triangle x_n)-f(x_n)]}{f(x_n+\triangle x_n)f(x_n)}\rv\\
&=O\lr \frac {|g|}f(x_n)\frac{f(x_n)}{f(x_n+\triangle x_n)}\lv\frac{\triangle x_n}{x_n}\rv^{1/2}\rr\\
&=O\lr \frac {|g|}f(x_n)\frac{f(x_n)}{f(x_n)+[f(x_n+\triangle x_n)-f(x_n)]}\lv\frac{\triangle x_n}{x_n}\rv^{1/2}\rr\\
&=O\lr \frac {|g|}f(x_n)\frac{1}{1+O\lr \lv\frac{\triangle x_n}{x_n}\rv^{1/2}\rr}\lv\frac{\triangle x_n}{x_n}\rv^{1/2}\rr=O\lr \frac {|g|}f(x_n)\lv\frac{\triangle x_n}{x_n}\rv^{1/2}\rr\\
&=O\lr \frac {|g|}f(x_n) \lr\frac{\log(\cdot)}{k_n}\rr^{1/4}\rr,
\end{split}
\end{equation*}
as $ \nty.$. The latter bound yields \eqref{psi_est1}. The corollary is proved. $\quad \square$
 
\section{Appendix}  
\label{appendix}
Let as before, $\na=\sharp \{i \ : X_i \le \xia ,\ i=1,\dots,n\}$, where $0<\alp <1$ is fixed. In this appendix we prove that
conditionally on $\na$ the order statistics $X_{1:n},\dots ,X_{\na:n}$ are distributed as
order statistics corresponding a sample of $\na$ i.i.d.~r.v.'s
with distribution function $F(x)/\alp$, $x\le \xia$. Though this fact is essentially known   (cf.,e.g.,~Theorem~12.4,~\cite{kallen}, cf.~also
\cite{gh2007}, \cite{hp}), we add a~brief proof of it. Let
$U_1,\dots,U_n$ be independent r.v.'s uniformly distributed on
$(0,1)$ and let $U_{1:n},\dots,U_{n:n}$ denote the corresponding
order statistics. Put $\nau=\sharp \{i \ : \ U_i \le \alp ,\ i=1,\dots,n\}$.
Since the joint distribution of the pair $\xin$, $\na$ is same as joint distribution of 
$ F^{-1}(\uin)$, $\nau\,$, it is enough to prove the assertion
for the uniform  
distribution.   
\begin{lemma}  
\label{Lemma_A} Conditionally given $\nau$, the order
statistics $U_{1,n},\dots,U_{\nau,n}$ are distributed as order
statistics corresponding to a sample of $\nau$ independent
$(0,\alp)$-uniform distributed r.v.'s.
\end{lemma}
\noindent{\bf Proof}.  $a).$ First consider the case $\nau= n$. Take
arbitrary $0<u_1\le \cdots \le u_n < \alp$ and write
\begin{equation*}
\label{ap_1}
\begin{split}
P(U_{1:n} \le u_1,\dots,U_{\nau:n}\le u_n \mid\nau =n)
=\frac{P(U_{1:n} \le u_1,\dots,U_{n:n}\le u_n)}{\alp^n}\\
=\frac{n!}{\alp^n}\int_0^{u_1}\int_{u_1}^{u_2}\dots
\int_{u_{n-1}}^{u_n} d\,x_1 d\,x_2\dots d\,x_n,
\end{split}
\end{equation*}
and the latter is $d.f.$ of the order statistics corresponding to
the sample of \ $n$ \ independent $(0,\alp)$-uniform distributed
r.v.'s. \ \  $b).$ Consider the case $\nau=k<n$. Let $F_{i,n}(u)=P(U_{i:n}\le u)$ be a~$df$ of $i$-th order statistic, put $P_n(k)=P(\nau=k)={n\choose k} p^k(1-p)^{n-k}$. Then we have
\begin{equation}
\label{ap_2}
P(U_{1:n} \le u_1,\dots,U_{\nau :n}\le u_k \mid \nau =k)
= \frac {P(U_{1:n} \le u_1,\dots,U_{k :n}\le u_k, U_{k+1:n}>p)}{P_n(k)}.
\end{equation}
The probability in the nominator on the r.h.s. of  \eqref{ap_2} is equal to
\begin{equation*}
\int_{\alp}^1 \!\!\! P\big(U_{1:n} \le
u_1,\dots,U_{k:n} \le u_k \mid U_{k+1:n}=v\big)\,d F_{k+1,n}(v),
\end{equation*}
and by the Markov property of order statistics the latter quantity equals
\begin{equation*}
\label{ap_3}
\begin{split}
&\phantom{ =}\int_{\alp}^1 \lr \frac{k!}{v^k} \int_0^{u_1}\int_{u_1}^{u_2}\dots \int_{u_{k-1}}^{u_k} d\,x_1
d\,x_2\dots d\,x_k\rr \,d F_{k+1,n}(v)\\
&= \frac{k!}{p^k}\lr \int_0^{u_1}\int_{u_1}^{u_2}\dots \int_{u_{k-1}}^{u_k} d\,x_1
d\,x_2\dots d\,x_k \rr\times p^k\int_{\alp}^1 \frac 1{v^k} \,d F_{k+1,n}(v),
\end{split}
\end{equation*}
and since $p^k\int_{\alp}^1 \frac 1{v^k} \,d F_{k+1,n}(v)=p^k\int_{\alp}^1 \frac{(1-v)^{n-k-1}}{B(k+1,n-k)} \,d v
= {n\choose k}p^k (1-p)^{n-k}=P_n(k)$, where $B(k+1,n-k)=k!(n-k-1)!/n!$, we obtain that conditional probability in  \eqref{ap_2} is equal
\begin{equation*}
\frac{k!}{\alp^k}\int_0^{u_1}\int_{u_1}^{u_2}\dots
\int_{u_{k-1}}^{u_k} d\,x_1 d\,x_2\dots d\,x_k,
\end{equation*}
which  corresponds to the $(0,\alp)$-uniform distribution. The lemma
is proved.$\quad \square$

\end{document}